\newif\ifdraft    
\documentclass[%
12pt,\ifdraft draft\fi%
]{article}
\usepackage{amsmath,amsfonts,amsthm,amssymb,amscd}
\usepackage{verbatim}	
\usepackage[matrix,arrow]{xy}
\CompileMatrices
\newsavebox{\diagrambox}
\newcommand\usediagram{\usebox{\diagrambox}}
  
\ifdraft
  
  \usepackage{showkeys}
\fi

\newcommand{\erdos}{Erd\H os}
\newcommand{\jonson}{J\'onsson}

\newcommand\bare{{\overline E}}
\newcommand\baru{{\overline U}}
\newcommand\barw{{\overline W}}
\newcommand\bcm{{\overline\cm}}
\newcommand\tcn{{\widetilde\cn}}
\newcommand\bct{{\overline\ct}}
\newcommand\cv{{\mathcal V}}
\newcommand\cw{{\mathcal W}}
\newcommand\tcu{{\widetilde\cu}}
\newcommand\wt{\widetilde}
\newcommand\cat{{^\frown}}
\newcommand\ppstree[3]{#2\oplus_{#1} #3} 
\let\pstree=\ppstree
\newcommand{\tlt}[1]{\lessdot^{#1}}
\newcommand\qi{quasi-iteration}
\newcommand\cmp{\,}		
\ifdraft
  \newcommand{\fref}[1]{((\ref{#1}, \pageref{#1}))}
  \newcommand{\xfootnote}[1]{ ******* \footnote{******* #1}}
\else
  \newcommand{\fref}[1]{}
  \newcommand{\xfootnote}[1]{}
\fi
\newtheorem{thm}{Theorem}[section]
\newtheorem{prop}[thm]{Proposition}
\newtheorem{lem}[thm]{Lemma}
\newtheorem{cor}[thm]{Corollary}
\newtheorem{claim}[thm]{Claim}
\theoremstyle{definition}
\newtheorem{defn}[thm]{Definition}
\theoremstyle{remark}

\newtheorem{question}{Question}
\newenvironment{case}[2]{\par\smallskip\noindent{\bf Case~#1 }({\sl #2\/}).\quad}{\par\smallskip}
\newcommand{\ga}{\alpha}
\newcommand{\gb}{\beta}
\newcommand{\gd}{\delta}
\renewcommand\gg{\gamma}
\newcommand\gga{\gamma}

\newcommand{\gk}{\kappa}
\newcommand{\gl}{\lambda}
\newcommand{\gs}{\sigma}
\newcommand{\gS}{\Sigma}
\newcommand{\gth}{\theta}
\newcommand{\gw}{\omega}
\newcommand{\ca}{{\mathcal A}}
\newcommand{\cb}{{\mathcal B}}
\newcommand{\ce}{{\mathcal E}}
\newcommand{\cf}{{\mathcal F}}

\newcommand{\cm}{{\mathcal M}}
\newcommand{\cn}{{\mathcal N}}

\newcommand{\cR}{{\mathcal R}}
\newcommand{\cs}{{\mathcal S}}
\newcommand{\ct}{{\mathcal T}}
\newcommand{\cu}{{\mathcal U}}
\newcommand{\mse}[1]{{\mathfrak{#1}}}
\newcommand\pcat{\oplus}	
\newcommand\union{\bigcup}
\newcommand\col{\operatorname{col}} 
\newcommand{\cof}{\operatorname{cf}}
\newcommand{\crit}{\operatorname{crit}}

\newcommand{\domain}{\operatorname{domain}}
\newcommand{\len}{\operatorname{len}}
\newcommand{\ot}{\operatorname{order\,type}}
\newcommand{\ord}{\operatorname{Ord}}

\newcommand{\range}{\operatorname{range}}
\newcommand{\ult}{\operatorname{ult}}
\newcommand{\card}[1]{\left|#1\right|}
\newcommand{\image}{\text{``}}
\newcommand{\ps}{\mathcal P}
\newcommand{\restrict}{\mathord\upharpoonright}
\newcommand{\set}[1]{\left\{\,#1\,\right\}}
\newcommand{\seq}[1]{\left\langle\,#1\,\right\rangle}

\begin{document}
\title     
{J\'onsson Cardinals, Erd\H os Cardinals,\\ and the Core Model} 
\author{W. J. Mitchell\thanks{This work was partially supported by grant number DMS-9306286
from the National Science Foundation.}}
\date{\today\ifdraft --- Preliminary\fi}
\maketitle
\begin{abstract}
We show that if there is no inner model with a Woodin cardinal and the Steel
core model $K$ exists, then every \jonson{} cardinal is Ramsey in
$K$, and 
 every $\gd$-\jonson{} cardinal is $\gd$-\erdos{} in $K$.

In the absence of the Steel core model $K$ we prove the same conclusion
for any model $L[\ce]$ such that either $V=L[\ce]$ is
the minimal model for a Woodin cardinal, or there is no inner model
with a Woodin cardinal and $V$ is a generic extension of $L[\ce]$.
\end{abstract}

\section{Introduction}
It is well known that every Ramsey cardinal is a \jonson{} cardinal
but surprisingly 
it appears to be unknown whether every \jonson{} cardinal is Ramsey.
It seems
unlikely that this converse implication holds in general, but Kunen
\cite{kunen.it-up} showed that it is true in $L[\mu]$, and this was
extended in \cite{mitchell.ramsey} (see also
\cite{jensen.applications}) to show 
that if $0^\dagger$ does not exist then every cardinal which is 
\jonson, even in $V$, is Ramsey in $K$.
In this note 
we  catch up with recent advances in core model theory by extending
this result to Steel's core model 
\cite{steel.core-model}, provided that this core model exists and
there is no model with a Woodin cardinal.  

Our main motivation in thinking about this problem came from an
interest in the possibility of a covering lemma for a model with a Woodin
cardinal.  This interest led to 
two extensions of the basic result:
First, we consider cardinals $\gk$ which are $\gd$-\jonson{} for 
a regular cardinal $\gd\le\gk$.  A $\gd$-\jonson{} cardinal, which is the 
same as a \jonson{} cardinal except that the the elementary
substructure is only required to have order type $\gd$, 
comes up in stationary tower
forcing~\cite{ma-so-woodin.determinacy} at a Woodin cardinal: 
$\ga$ is $\gd$-\jonson{} if and only if there is 
a condition in the stationary tower forcing which forces that
$i(\gd)=\ga$, where $i$ is the generic 
embedding.  We show, under the same conditions as above, that
every $\gd$-\jonson{} cardinal is $\gd$-\erdos{} in $K$.  To some extent
we are following Jensen in this:
in
\cite{jensen.applications}, Jensen extended \cite{mitchell.ramsey} by
showing that if there are no models with a measurable cardinal then
every  $\gd$-\erdos{} cardinal is $\gd$-\erdos{} 
in the Dodd-Jensen core model for one measurable cardinal.
 
The relevance of these results to Woodin cardinals would seem to be
limited by the fact that the relevant Steel core models do not exist
in the presence of a Woodin cardinal, and (so far as we know) may not
exist even if there is no model with a Woodin cardinal.
However, if $V=L[\ce]$ or $V$ is a generic extension of a model
$L[\ce]$ then the model $L[\ce]$ looks like a core model. We use this
resemblance to show that if $L[\ce]$ is a minimal model for a Woodin
cardinal then $L[\ce]$ satisfies that every $\gd$-\jonson\
cardinal is $\gd$-\erdos, and that if $L[\ce]$ does not contain a
class model for a 
Woodin cardinal then the same is true of cardinals which are
$\gd$-\jonson{} in a generic extension of $L[\ce]$.

Before stating the main theorem, we give precise definitions of
$\gd$-\jonson{} and $\gd$-\erdos{} cardinals.

\begin{defn}
If $\gd$ and $\gk$ are cardinals with $\gd<\gk$, then $\gk$ is said
to be {\em
$\gd$-\jonson\/} if for each first order structure $\ca$ in a
countable language with universe $\gk$  there is an elementary
substructure $\ca'\prec\ca$ with universe $A'$ such that
$\ot(A')=\gd$. 

We say that $\gk$ is $\gk$-\jonson{} if it is \jonson, that is, if
for every structure $\ca$ as above there is an elementary substructure
$\ca'\prec\ca$ with universe $A'$ such that $\card{A'}=\gk$ but
$A'\not=\gk$.
\end{defn}

The following definition is due to Baumgartner \cite{baumgartner.erdos}:
\begin{defn}
If $\gd$ and $\gk$ are cardinals, with $\gd\le\gk$, then $\gk$ is {\em
$\gd$-\erdos\/} if for every structure $\ca$ in a countable language
with universe~$\gk$ and for every closed and unbounded subset $C$ of
$\gk$ there is a set $D\subset C$ of order type $\gd$ which is a
\emph{normal set of indiscernibles} for $\ca$.
\end{defn}

By a \emph{normal set of indiscernibles} we mean a set $D$ such that 
for every $n$-ary function $f$ which is definable in $\ca$ without
parameters,   
either $f(d_0,\dots,d_{n-1})\ge d_0$ for every 
$\vec d=\seq{d_0,\dots,d_{n-1}}\in[D]^{n}$ or else      
the value of $f(\vec d)$ is constant for
$\vec d\in [D]^{n}$. 
It should be noted this is equivalent to the seemingly stronger
definition which requires the same property to hold for 
all functions which are definable in $\ca$ from
parameters smaller than $d_0$.

It is shown in \cite{baumgartner.erdos} that every Ramsey cardinal
$\gk$  is $\gk$-\erdos.

\begin{thm}[Main Theorem]
Suppose that $L[\ce]$ is a fully iterable model constructed from an good
sequence of extenders as in \cite{mitchell-steel.inner-model}, and
assume that any of the following three conditions are true:
\begin{enumerate}
\item There is no class model with a Woodin cardinal, and $L[\ce]$ is
equal to the Steel core model $K$.
\item There is no class model with a Woodin cardinal, and $V$ is a generic extension of $L[\ce]$.
\item $V=L[\ce]$ is a minimal model for a Woodin cardinal: that is,
if 
$\eta=\len(\ce)$ then $\eta$ is Woodin in $L[\ce]$, but all of the 
models $L[\ce\restrict\ga]$ for $\ga<\eta$ are fully iterable.
\end{enumerate}

Then every cardinal $\gk$ which is a $\gd$-\jonson{} cardinal in $V$,
where $\gd$ is any 
uncountable regular
cardinal, is $\gd$-\erdos{} in $L[\ce]$.
\label{t:j->erdos}
\end{thm}

The statement in clause~(3) that a model $L[\ce]$ is minimal requires
not only that there is no $\ga<\len(\ce)$ such that 
$L[\ce\restrict\ga]$ satisfies that
$\ga$ is Woodin, but also that there is no iteration tree $\ct$ on
$L_{\ga}[\ce\restrict\ga]$ such that $L[\ce^{\ct}]$ satisfies that
there is a Woodin cardinal, where $\ce^{\ct}$ is the extender sequence
stabilized by $\ct$.  This means that $\ce^{\ct}$ is the union of the
set of sequences $\ce_{\nu}\restrict\rho_\nu$, for $\nu$ on the main
branch of $\ct$, where $\ce_{\nu}$ is the extender sequence of the
$\nu$th tree of $\ct$ and
$\ce_{\nu}\restrict\rho_\nu=\ce_{\nu'}\restrict\rho_\nu$ for all
$\nu'>\nu$ in the main branch of $\ct$.
The assumption  in hypotheses~(1) and~(2) that there is no class
model with a Woodin 
cardinal is required 
to ensure iterability: we need to know that every iteration tree in $V$
on $L[\ce]$ has a branch in $V$. 

We begin by reducing hypothesis~(3) of the theorem to hypothesis~(2), thereby
eliminating the one case in which there is a Woodin cardinal.  

\begin{lem}
\label{l.reduce}
Suppose that $L[\ce]$ is a minimal model for a Woodin cardinal and
that $\gk$ is $\gd$-\jonson{} but not $\gd$-\erdos{} in $L[\ce]$ for some
$\gd\le\gk$ which is regular and uncountable in $L[\ce]$.
Then there is $\ga<\len(\ce)$ so that
$L[\ce\restrict\ga]$ also satisfies 
that $\gk$ is $\gd$-\jonson{} but not $\gd$-\erdos.
\end{lem}
\begin{proof}
Let $L[\ce]$, $\gk$ and $\gd$ be as in the hypothesis of the lemma,
and let $\eta=\len(\ce)$.  Then $\gk\le\eta$, since there are no
$\gd$-\jonson{} cardinals in $L[\ce]$ above~$\eta$.

If $\gk<\eta$ then  we can take $\ga=\gk^{+}<\eta$.  Then 
$\gk$ is $\gd$-\jonson, but not $\gd$-\erdos, in $L[\ce\restrict\ga]$ since
$\ps(\ga)\cap L[\ce]=\ps(\ga)\cap L[\ce\restrict\ga]$.

If $\gd<\gk=\eta$ then the hypothesis of the lemma is false since 
 the set
of  measurable cardinals is  
stationary in $\eta$, which implies that $\eta$ is $\gd$-\erdos{}
 in~$L[\ce]$.  

We will finish the proof of the lemma by showing that if
$\gd=\gk=\eta$ then  the structure
$L_{\eta}[\ce]$ is a \jonson{} algebra in $L[\ce]$, so that 
the hypothesis of the lemma is false since 
$\eta$ is not \jonson.

\begin{lem}\footnote{NOTE: This replaces an earlier incorrect proof.
It is so easy that it 
must have been known before, and I  will probably replace this by a
reference, once I find one.}
Suppose $V=L[A]$ where $A\subset\gk$ and $\gk$ is a regular cardinal.
Then $L_{\gk}[A]$ is a \jonson\ algebra, and hence  $\gk$ is not a
\jonson\ cardinal. 
\end{lem}
\begin{proof}
Suppose to the contrary that there is an $X\prec L_{\gk}[A]$ such that
$\card X=\gk$ 
but $\gk\not\subset X$.  Choose such an $X$ with
$\ga=\inf(\gk\setminus X)$ as small as possible, and let $\pi\colon
L_{\gk}[A']\cong X\prec L_{\gk}[A]$.   

We claim that $\pi$ can be extended to an elementary embedding
$\tilde\pi\colon L[A']\to \ult(L[A'],\pi,\gk)=L[A]$.  For this it is
sufficient to show that $\ult(L[A'],\pi,\gk)$ is defined and well
founded.   

To show that  $\ult(L[A'],\pi,\gk)$ is defined we need to show that
every bounded subset of $\gk$ in $L[A']$ is a member of $L_{\gk}[A']$.
Suppose $x\subset\gd<\gk$ and $x\in L_{\xi}[A']$.  Then there is
$Z\prec L_{\xi}[A']$ with $x\in Z$ such that $\card Z<\gk$ and
$\gd\le\eta=Z\cap 
\gk=\sup(Z\cap\gk)$.   Let $k\colon L_{\xi'}[A'\cap\eta]\cong Z$ be
the inverse of the transitive collapse.  Then 
$x\in L_{\xi'}[A'\cap\eta]\in L_{\gk}[A']$.

A similar argument shows that $\ult(L[A'],\pi,\gk)$ is well founded: 
Otherwise pick functions $f_n$ and ordinals $\ga_n<\gk$ so that 
$\seq{\tilde\pi(f_n)(\ga_n):n<\gw}$ is a descending sequence  of
ordinals.
Pick $Z$ as above with
$\set{\ga_n:n<\gw}\cup\set{f_n:n<\gw}\subset Z$, and define $k\colon
L_{\xi'}[A'\cap\eta]\cong Z$ as before.  Then $L_{\xi'}[A'\cap\eta]\in
L_{\gk}[A']$, and $\pi\bigl(L_{\xi'}[A'\cap\eta]\bigr)$ is ill
founded, which is absurd since it is a member of $L_{\gk}[A]$.

Since $k(\ga)>\ga$ and $\tilde\pi\colon L[A'] \to L[A]$ is an
elementary embedding, there is a set $X'\in L[A']$ such that 
$\card X'=\gk$ and $X'\prec L_{\gk}[A']$ but $\ga'=\inf(\gk\setminus
X')<\ga$.
Now set $X''=\pi\image X'$.  Then
$\card{X''}=\card{X'}=\gk$ and $X''\prec L_{\gk}[A]$, but $\ga'\notin
X''$.  This contradicts the minimality of $\ga$, and hence completes
the proof of the lemma.
\end{proof}
This finishes the proof of lemma~\ref{l.reduce}.
\end{proof}

It follows that it is sufficient to prove that the conclusion of
the main theorem follows from hypothesis~(1)
or~(2) of that theorem.  Before doing so we complete the
introduction by discussing some of the notation and background theory
used in the proof.

\subsubsection*{Iteration trees, $L[\ce]$, and the Core model}

The basic sources are \cite{mitchell-steel.inner-model} for $L[\ce]$
and \cite{steel.core-model}  for Steel's core model $K$.  The primary aim
of this section is to clarify the notation which we will use and to
state some basic results concerning iteration trees, 
$\varphi$-minimal structures $L[\ce]$ and
the core  $K$.

 A \emph{phalanx} of length $\gth+1$ is a pair
$(\vec\rho,\vec\cm)$
where $\rho=(\rho_\nu:\nu<\gth)$ is a continuous increasing sequence of
ordinals and $\vec\cm=(\cm_{\nu}:\nu\le\gth)$ is a sequence of premice
such 
that if $\ga<\gb\le\gth$ then $\cm_{\ga}$ agrees with $\cm_{\gb}$ up
to~$\rho_{\ga}$. 

We will write $\cu\oplus\cu'$ for a phalanx which has $\cu$ as an
initial segment.  For the special case where $\cu'$ is a phalanx of
length~1, that is, a single premouse $\cR$, we will write
$\cu\oplus_{\gl}\cR$ for the phalanx obtained by truncating $\cu$ as
necessary, and then adding $\cR$ onto the end.  If
$\cu=(\vec\rho,\vec\cm)$ then this means that 
\[
\cu\oplus_{\gl}\cR=\left(\vec\rho\restrict\nu\cat\gl,\vec\cm\restrict(\nu+1)\cat\cR\right)
\]
where $\nu$ is the least ordinal such that $\rho_{\nu}>\gl$.

Any iteration tree has an underlying phalanx, which is simply the
sequence of models $\cm_{\nu}$ and ordinals $\rho_{\nu}$ of the tree.
We will normally use the same symbol for a tree and for its associated
phalanx.   

We will use the term \emph{iteration tree} both for a normal iteration
tree, with a single root, and for an \emph{iteration tree on a phalanx},
which Steel calls a \emph{pseudo-iteration tree}.
An iteration tree $\ct$ on a phalanx
$\cu=(\vec\rho,\vec\cm)$ is an iteration tree in the
normal sense, except that the underlying phalanx of $\ct$ has the form
$\cu\oplus\ct'$ and the roots of $\ct$ are exactly the members of
$\cu$.

Thus, suppose that
$F_\nu$ is the extender to be used at stage $\nu<\len(\ct)$ of
the construction of $\ct$, 
and  that 
$\eta=\crit(F_{\nu})<\rho_\ga$ for some $\ga<\len(\cu)$.  Let $\ga$ be
the largest ordinal such that $\rho_{\ga}\ge\eta$.  Then 
$\nu+1$
is at the second level of the tree, being  an immediate successor to
the root $\ga$, and $N_{\nu+1}=\ult(\cm_{\ga}, F_\nu)$.

We will use the symbol $\tlt{\ct}$ to denote the tree ordering on an
iteration tree $\ct$; regarding it as an ordering on either the models
of the tree or on their indices depending on which is more convenient
for the exposition.  Thus if $\cm_{\nu} $ is the $\nu$th model of
the tree $\ct$ then the two formulas $\nu\tlt{\ct}\nu'$ and
$\cm_{\nu}\tlt{\ct}\cm_{\nu'}$ mean the same thing. 

\paragraph{Iteration trees and comparisons.}
If $\ct$ is an iteration tree without drops on its main branch then
we write $i^{\ct}$ for the embedding along the main branch.  
If $\ct$ is a normal iteration tree and $P$ is its last model 
then $i^{\ct}\colon\cn_0\to P$.  If $\ct$ is an iteration tree on a
phalanx $\cu$ then 
$i^{\ct}\colon\cn_{\nu}\to P$, where $\cn_\nu$ is the unique member of $\cu$
such that $\cn_{\nu}\tlt{\ct} P$.   

If $\cu$ and $\cu'$ are two phalanxes then the  comparison of $\cu$ and
$\cu'$ yields iteration trees $\ct$ on $\cu$ and $\ct'$ on $\cu'$ such
that if $P$ and $P'$ are the last models of $\ct$ and $\ct'$ then one
of the models $P$ and $P'$ is an initial segment of the other. 

The following lemma gives some basic standard facts about this
comparison.  See, for example, \cite{mitchell-steel.inner-model}
\begin{lem}
\label{l.compare}
Suppose that $\cu$ and $\cu'$ are compared using trees $\ct$ and
$\ct'$, with last models $P$ and $P'$.
\begin{enumerate}
\item
At most one of the trees $\ct$ and $\ct'$ has a drop on its main
branch, and  if $\ct$ has such a drop then $P'$ is an initial segment
of $P$.
\item
\label{l.comp3}  
Suppose that the trees $\ct$ and $\ct'$ are not both trivial and that 
there is no drop in the main branch of either tree, and  let
$\ga=\min(\crit(i^{\ct}),\crit(i^{\ct'}))$.  Then   
$i^{\ct}\restrict\ps(\ga)\not= i^{\ct'}\restrict\ps(\ga')$.   
\end{enumerate} 
\end{lem} 

The proof of the conclusion of the main theorem is essentially the
same under the assumption of either hypothesis~(1) or~(2).  The next two
lemmas express the facts we need about the relevant models:
The notion of $\varphi$-minimality and lemma~\ref{l.phi-minimal} are used for
hypothesis~(1), while lemma~\ref{l.steel-core-model} is used for hypothesis~(2).

\begin{defn} If $\varphi$ is a sentence of set theory then a model
$L[\ce]$ is said to be \emph{$\varphi$-minimal} if
$L[\ce]\models\varphi$, but no model $L[\ce\restrict\ga]$ or
$L_{\gb}[\ce\restrict\ga]$ 
with $\ga<\len(\ce)$ satisfies $\varphi$.
\end{defn}

Notice that
only the first two clauses of the following lemma assume
$\varphi$-minimality.

\begin{lem}  Suppose that $\ce$ is a good sequence
of extenders as in \cite{mitchell-steel.inner-model} and $L[\ce]$ is fully
iterable.
\label{l.phi-minimal}
\begin{enumerate}
\item If $i\colon L[\ce]\to P$ is elementary and $L[\ce]$ is
$\varphi$-minimal for some sentence $\varphi$ 
then, $P$ is also $\varphi$-minimal.
\item Suppose that  $L[\ce]$ is $\varphi$-minimal and $\ct$ is an
iteration tree with last model $P$.  If there are no drops on the main
branch of $\ct$ then $P$ is $\varphi$-minimal; while if there are any drops
on the main branch of $\ct$ then neither $P$  nor any initial segment
of $P$ satisfies $\varphi$.
\item Suppose that $i\colon L[\ce]\to P$ and $j\colon L[\ce]\to Q$ are
elementary embeddings 
which are definable in a generic extension of $L[\ce]$, and that $P$
is an initial segment of~$Q$.  Then $P=Q$ and $i=j$.
\item More generally, suppose that the elementary embeddings 
\[ 
L[\ce]\mathop{\longrightarrow}\limits^k
L[\ce']
\stackrel{\displaystyle\mathop{\longrightarrow}\limits^i}
         {\mathop{\longrightarrow}\limits_j}P
\]
 are  definable in a generic
extension of $L[\ce]$ and that $k$ is generated by
$\rho=\min(\crit(i),\crit(j))$, that is, 
\[
L[\ce']=\set{k(f)(\nu):f\in
L[\ce]\text { and } \nu<\rho}.
\]  Then $i=j$.
\item
Suppose that $\cm$ is a mouse with projectum $\ga$ which agrees with
$L[\ce]$ up to $\ga$, and that the phalanx $(L[\ce],(\ga,\cm))$ is
fully iterable.  If $\cm$ is a member of a generic extension of
$L[\ce]$ then $\cm$ is a member of $L[\ce]$.
\end{enumerate}
\end{lem}
\begin{proof}
The proof of clause~(1) is immediate, and  clause~(2) can be proved by a
straightforward induction on the length of the iteration tree $\ct$.

Suppose for the sake of a contradiction that clause~(3) is false.
The assertion that clause~(3) is false is a first order statement
$\varphi$ over $L[\ce]$, so we can suppose that $L[\ce]$ is
$\varphi$-minimal. 
Now pick a partial order $\mathbb P\in L[\ce]$, a $L[\ce]$-generic
set $G\subset \mathbb P$ and embeddings $i\colon L[\ce]\to P$ and
$j\colon L[\ce]\to Q$ in $L[\ce][G]$ witnessing
the failure of clause~(3).  Let $x$ be the least set, in the order
of construction of $L[\ce]$, such that $i(x)\not=j(x)$.  We may
suppose that $\mathbb P$, $G$ and $i$ and $j$ were chosen so that $x$
is as small as possible; thus $x$ is definable in $L[\ce]$.

Then $P$ and $Q$
are both 
$\varphi$-minimal by clause~(1), and hence $P=Q$. 
But then $i(x)=j(x)$ since $x$ is definable in $L[\ce]$. This  contradicts the
choice of $x$ and hence completes the proof of clause~(3).

Clause~(4) follows from clause~(3): every member $x$ of $L[\ce']$ can be
written in the form $x=k(f)(a)$ where $f\in L[\ce]$ and
$a\in[\ga]^{<\gw}$.  Then 
\begin{equation*} 
i(x)=i(k(f)(a))=(i\cmp k(f))(i(a))=      
(j\cmp k(f))(j(a))=j(k(f)(a))=j(x), 
\end{equation*} 
since $i\cmp k=j\cmp k$ by clause~(3) and $i(a)=j(a)=a$.

 We can use standard arguments to prove clause~(5) from
clauses~(1--4).  We will give a fairly complete proof here in order to
remind the reader of the techniques which will be applied later in the
paper in slightly different contexts.  As the following diagram
indicates, we begin by comparing $L[\ce]$ with the phalanx
$(\ga,(L[\ce],\cm))$, using trees $\ct$ on $L[\ce]$ and $\cu$ on
$(\ga,(L[\ce],\cm))$.  We use wavy arrows in the diagram 
since we do not know whether the
indicated embeddings exist:
\begin{lrbox}{\diagrambox}
\xymatrixrowsep={1ex}
\xymatrix{
*!R{L[\ce]}          \ar@{~>}[r]^{\ct}     & P\\
\\
*!R{(\ga,(L[\ce],\cm))} \ar@{~>} [r]^{\cu}  &Q
}
\end{lrbox}
\begin{equation*}
\usediagram
\end{equation*}
Thus the first two models of the tree $\cu$ are
$\cn_0=L[\ce]$ and $\cn_1=\cm$.  The tree $\cu$ has two roots $0$,
and $1$, which means that  
 $0$ and $1$ are incomparable in the tree
ordering  $\tlt{\cu}$ of $\cu$, while
for every $\nu<\len(\cu)$ either $0\tlt{\cu}\nu$ or $1\tlt{\cu}\nu$.
We will  
say for short that every model $\cm_{\nu}$ of $\cu$ is either above
$L[\ce]$ or above $\cm$.

The ordinal $\ga$ is used as if
it were the length of an extender used to obtain $\cn_1$.  This means
that if $\nu<\len(\cu)$ then $0$ is an immediate
predecessor of $\nu+1$ in $\cu$, so that 
$N_{\nu+1}=\ult(L[\ce],F_{\nu})$,  
 if and only if $\crit(F_\nu)<\ga$,
where $F_\nu\in\cn_\nu$ is  the extender which is to be used to
define $\cn_{\nu+1}$.

First, notice that if the last model $Q$ of $\cu$ lies above $\cm$
and is an initial segment of $P$ then $\cm\in L[\ce]$, as required.
To see 
this, first note that since $Q$ is an
initial segment of $P$, lemma~\ref{l.compare}(1)  implies that there is 
no drop on the 
main branch of $\cu$, so that the embedding $i^{\cu}$ is defined.
Furthermore $i^{\cu}\restrict\ga$ is the identity 
since $Q$ is above $(\ga,\cm)$ in $\cu$.  
Hence the master code $A$ of $\cm$ is
still definable in $Q$, and hence is a member of $P$.  But since
$L[\ce]$ and $\cm$ agree up to $\ga$, it follows that $A\in L[\ce]$,
and thus $\cm\in L[\ce]$ since it is definable from $A$.

Thus it is sufficient to show that $Q$ lies above $\cm$ in $\cu$, and
that $Q$ is an initial segment of $P$.
We will first show that $Q$ lies above $\cm$.  Suppose to the contrary
that $Q$ lies above $L[\ce]$.  Then there is no drop in the main
branch of either tree: suppose for example that there is a drop in the
main branch 
of $\ct$.  Then no initial segment of $P$ satisfies $\varphi$.  But $Q$ is
an initial segment of $P$ and hence also fails to satisfy $\varphi$.
This
is a contradiction since the embedding $i^{\cu}\colon L[\ce]\to Q$ is
defined and so $Q\models\varphi$ by elementarity.   Thus the main branch of
$\ct$ does not drop, and the same argument shows that the main branch
of $\cu$ does not drop either.  It follows by clause~(3) that $i=j$, but this
contradicts clause~\ref{l.comp3} of lemma~\ref{l.compare}.  Thus $Q$ lies above
$\cm$ in $\cu$.

Now we can finish the proof by showing that $Q$ is an initial segment
of $P$.   Suppose to the contrary that $P$ is a proper initial segment
of $Q$.  This implies that there are no drops on the main branch of
$\ct$, so $P$ is a proper class.  Then $Q$ is a proper class, and since
it  lies above $\cm$, which is a set,
it follows that $\len(\cu)=\ord$ and that there is a closed
unbounded class $C$ of ordinals which are critical points of the
embeddings along the main branch of $\cu$, and hence are inaccessible
in $Q$.  Furthermore 
there is a closed unbounded subclass $C'\subset C$ such that
$i^{\ct}\image\ga\subset\ga$ for $\ga\in C'$, and   since the universe is a
set generic extension of $L[\ce]$ there is a closed unbounded subclass
$C''\subset C'$ which is definable in $L[\ce]$.  If $\ga$ is the
$\gw$th member of $C''$ then $\cof(\ga)=\gw$  in $L[\ce]$, and hence
in $P$.  This is a contradiction because $\cof^{Q}(\ga)=\ga$, but $Q$
contains $P$.

This completes the proof of lemma~\ref{l.phi-minimal}.
\end{proof}

\medspace

The next lemma is used instead of lemma~\ref{l.phi-minimal} 
to prove the conclusion of the main theorem from hypothesis~(2).  In
Steel's terminology, it asserts that the initial 
segments of the core model are \emph{very sound}.
\begin{lem}[Steel \cite{steel.core-model}]
Let $K=L[\cf]$ be the Steel core model and let $\gl$ be any
ordinal.  Then there is a model $W=L[\ce]$ such that
\label{l.steel-core-model}
\begin{enumerate}
\item 
$\ce\restrict\gl = \cf\restrict\gl$.
\item ($W$ is universal)
If $W$ is compared with any iterable phalanx $\cu$, then the last
model of the tree on $\cu$ is an initial segment of the last model of
the tree on $W$.  Furthermore, if $i\colon W\to W'$ is elementary and
$W'$ is iterable then $W'$ is also universal.
\item 
In particular, if  $\cm$ is a mouse with projectum $\ga\le\gl$ which
agrees with
$W$ up to $\ga$, and if the phalanx
$(\ga,(W,\cm))$ is fully iterable,  then $\cm$ is an initial
segment of $W=L[\ce]$, that is, $\cm=J_{\nu}[\ce]$ for some
ordinal~$\nu$.  
\item
If $i,j\colon W\to P$ for some iterable model $P$ then $i\restrict
K_{\gl}=j\restrict K_{\gl}$.
\end{enumerate}
\end{lem}

\section{Proof of the main theorem}
\subsection*{Notation and Summary}
As pointed out earlier, it is sufficient to prove the conclusion of
theorem~\ref{t:j->erdos} from hypotheses~(1) and~(2).  For the rest of
this paper we assume 
that $W=L[\ce]$ is an iterable model, that there is no iterable class
model with a Woodin cardinal, and that $\gd\le\gk$ are cardinals in
$W$ such that $\gd$ is regular and uncountable and $\gk$ is
$\gd$-\jonson.  Furthermore we assume for hypothesis~(1) that  $V$ is a generic
extension of $W$, so that we can use
lemma~\ref{l.phi-minimal}, and that $W$ is  $\varphi$-minimal for the assertion that the conclusion of the main
theorem fails in some generic extension of $W$.
 We assume for hypothesis~(2) 
that the core model $K$ exists, and that $W$ is a model agreeing with
$K$ up to $\gk^+$ which 
satisfies the conditions of lemma~\ref{l.steel-core-model}.
In either case we
will show that $\gk$ is $\gd$-\erdos{} in $W$.  This gives a direct
proof of the conclusion of the theorem from hypothesis~(2), and a proof by
contradiction from hypothesis~(1). 

Since we are trying to show that $\gk$ is $\gd$-\erdos{} in $L[\ce]$,
let us fix an arbitrary structure  $\ca\in L[\ce]$ in a countable language 
with universe $\gk$, and let $C\in
L[\ce]$ be a closed unbounded subset of $\gk$.
We will
find a set $D\in L[\ce]$ of  normal indiscernibles for $\ca$ such that
$D\subset C$ and $\card D=\gd$.  
Set $\gl=(2^{\gk})^{+}$.

\begin{prop}
\label{Xdef}
There is a set $X$ satisfying  $\set{\gd,\gk,C,\ca}\subset X$ and
$(X,\ce)\prec(H_\gl ,\ce)$ such that $\gd\not\subset X$ but
$X\cap\gk$ has order-type $\gd$.
\end{prop}
\begin{proof}
First take an elementary substructure $(X^*,\ce)$ of $(H_\gl ,\ce)$
such that
$\set{\gd,\gk,C,\ca}\cup\gk\subset X^*$ and $\card{X^*}=\gk$.  
Now let
$f\colon\gk\cong X^*$ and use $f$ to code $(X^*,\ce)$ into a structure
$\cb$ with
universe $\gk$ such that if $Z\subset\gk$ is the universe
of an elementary substructure of $\cb$ then $f\image Z$ is the
universe of an
elementary substructure of $(X^*,\ce)$ with the property that
$\gk \cap f\image 
Z=Z$ and $\set{\gd,\gk,C,\ca}\subset f\image Z$.

  If $\gd=\gk$ then $\gk$ is \jonson{} and hence $\cb^*$ has a
elementary substructure with universe $Z$ such that $\card Z=\gk$ but
$Z\not=\gk$.   If  $\gd<\gk$ then $\gk$ is $\gd$-\jonson\
and there is an elementary substructure $\cb$ with universe $Z\not\subset\gd$ such that 
$\ot{Z}=\gd$, so that $\ot{Z\cap\gd}<\gd$ and hence $\gd\not\subset\cb^*$.
In either case $X=f\image Z$  is the universe of 
an elementary substructure of $(X^*,\ce)$
with the required properties.
\end{proof}

Now let $X$ be as given by the proposition and let $N$ be a transitive
set with $\pi\colon N\cong X$, so that $\crit(\pi)<\gd$.
This situation is similar to the situation at the start of the proof
of the 
weak covering lemma 
\cite{mit-sch-ste.CovLemWoodin}, and it will be useful to compare the
two proofs.  Our substructure $X$ differs from 
that used in the proof of the covering lemma in two significant ways.
The first is that we cannot assume that 
$^{\gw}X\subset X$, 
as in the proof of the covering lemma
(nor have we been able to use Fodor's lemma to
avoid countable closure, as Dodd and Jensen do) 
The other difference partially counterbalances the first:
$\gd=\pi^{-1}(\gk)$ is  a regular cardinal.
The closure condition $^{\gw}X\subset X$ 
is used several times in the proof of the covering
lemma, and we will deviate from the proof of the covering lemma only
when it is necessary to work around this lack of closure. 

As in the proof of the covering lemma,  we set
$\barw=\pi^{-1}(W)$, and 
compare the two models  $W$ and  $\barw$ using iteration
trees $\ct$  on $W$ and  $\cu$ 
on $\barw$, continuing this comparison  until the final models of
the two trees agree  up to  $\pi^{-1}(\gk)$.
At this point the proof of the weak covering lemma uses a
rather complicated induction to reach two important conclusions:
(i)~the sequence $\barw$ is never moved in the comparison, so that the tree $\cu$ is
trivial, and 
(ii)~if $\cm_{\phi}$ is the final model of $\ct$
then the model $\cR_{\phi}=\ult(\cm_{\phi},\pi,\gk)$ is iterable.  

We do not reach either of these two conclusions. The inability to
prove that $\cu$ is trivial is merely a nuisance; it will be dealt
with in the proof but for clarity we ignore it in this summary.  Our
inability to prove that $\cR_\phi$ is iterable, on the other hand, 
requires a fundamental change in the proof.  To see what  changes
are necessary, let us look at the two basic cases which come up in the
proof of the covering lemma:  

\begin{case}
{1}{$\cm_{\phi}$ is a set}
In this case there must be a drop somewhere along
the main branch $b$ of $\ct$, so that there is a $\nu<\phi$ in $b$
such that the
$\nu$th model $\cm_{\nu}$ of $\ct$ has cardinality less than~$\gd$.
In this case we can use the fact that $\pi^{-1}(\gk)=\gd$,  a
regular cardinal, to 
show  that $\ct$ has length $\gd$ and hence generates the required set
of indiscernibles of order type $\gd$.
There is no need for $\cR_{\phi}$ to be iterable in this case.
\end{case}

\begin{case}
{2}{$\cm_{\phi}$ is a weasel, that is, a proper class}
In this case we
will use an argument taken from the proof of the weak covering
lemma in \cite{mit-sch-ste.CovLemWoodin} 
to show that 
$i^{\ct}_{0,\phi}(\rho)>\gd$, where $i^{\ct}_{0,\phi}$ the embedding  along the main
branch of $\ct$  and $\rho=\crit(i^{\ct}_{0,\phi})$.  Since $\ct$ only
uses extenders of length less than $\gd$ we can again use the fact
that $\gd$ is a regular cardinal to conclude that $\ct$ has length
$\gd$ and hence generates the desired set of indiscernibles.  The
problem is that the argument taken from the covering lemma depends
heavily on the assumption that 
$\cR_{\phi}=\ult(\cm_{\phi},\pi,\gk)$ is iterable. We will work around
this difficulty by noticing that if $\cR_{\phi}$ is not iterable then
this failure must have been evidenced in some earlier model $\cm_\nu$
on the main branch of the tree $\ct$, and in fact 
in a structure  $Q\prec\cm_{\nu}$ with $\card Q<\gd$.
In this case we will modify the construction of the
tree $\ct$ by dropping  at stage $\nu$ to the premouse $Q$.

If there is any drop on
the main branch of the modified tree 
$\ct$ then we are in case~(1) and there is no need
for $\cR_{\phi}$ to be iterable.  On the other hand, if
there is no such drop then $\cR_{\phi}$ is iterable
and we can use the argument from the covering lemma.

For lack of a better term will will call the modified tree $\ct$ a
\emph{\qi{} tree}.
\end{case}

\medskip
We are now almost ready to begin the first half of the actual proof,
which is the 
construction of the trees $\ct$ on $W$ and $\cu$ on $\barw$.
 We also use the embedding $\pi$ to copy $\cu$ to a tree
$\tcu$ on $W$, using the  shift lemma of Martin and Steel
\cite{martin-steel.iteration-trees}. 
The trees $\cu$ and $\tcu$ are ordinary iteration trees in the sense of
\cite{mitchell-steel.inner-model} and hence  present  no difficulties
concerning iterability; however
 the tree $\ct$ is not a standard iteration tree and hence 
requires special treatment.
The verification of the iterability of $\ct$ is in 
lemma~\ref{l.strategy},  which is one of two lemmas which are
proved after the description of the construction of the trees.

The following diagram gives the maps between the trees $\ct$, $\cu$
and $\tcu$:
\begin{lrbox}{\diagrambox}
\xymatrix@R=.5pc{
*+!/r.7em/{\ct\colon\cm_{\nu}}  \ar@{~>}[r]^{i_{\nu,\nu'}}
                 &{\cm_{\nu'}}\\
*+!/r.7em/{\cu\colon\cn_{\nu}}  \ar@{~>}[r]^{j_{\nu,\nu'}}
                            \ar[dd]_{\pi_{\nu}}
                 &{\cn_{\nu'}} \ar[dd]_{\pi_{\nu'}}\\ 
						  \\
*+!/r.7em/{\tcu\colon\tcn_{\nu}} \ar@{~>}[r]^{\tilde\jmath_{\nu,\nu'}}
                 &{\tcn_{\nu'}}
}
\end{lrbox}
\begin{equation*}
\usediagram
\end{equation*}
The horizontal maps are only defined if $\nu'$ is above $\nu$ in the
relevant tree and there is no drop in the branch between~$\nu$ and~
$\nu'$.

\subsection*{The construction of $\ct$, $\cu$ and $\tcu$.}
For the rest of this proof $\pi\colon N\cong X\prec (H_{\gl},\ce)$
will be as in lemma~\ref{Xdef}.  We will write $\barw=\pi^{-1}[X\cap
W]$.

The trees $\ct,\cu$ and $\tcu$ are defined by recursion on their
lengths.  We are not using padded iteration trees, and hence the trees
need not have the same length.  We will write $\cm_{\nu}$ for the
$\nu$th model of $\ct$, and we will write $\cn_{\nu}$ and $\tcn_{\nu}$
for the $\nu$th model of $\cu$ and $\tcu$ respectively. 
The construction starts with $\cm_0=W$, $\cn_0=\barw$, $\wt\cn_0=W$,
and $\pi_0=\pi\colon \cn_0\to\wt\cn_0$. 

Suppose that
during the course of the recursion we have already defined an initial
segment $\ct\restrict\phi$ of $\ct$ and initial segments
$\cu\restrict\gth$ and $\tcu\restrict\gth$ of $\cu$ and $\tcu$
respectively.  At the next stage of the recursion we will extend one
or both of the trees $\ct$ and $\cu$. There are three cases:

\begin{case}{1}{At least one of $\phi$ or $\gth$ is a limit ordinal}
  If $\gth$ is a limit ordinal then we extend the tree $\cu\restrict\gth$ 
by taking the unique well founded branch $b$ of
$\cu\restrict\gth$.  This
unique well founded branch exists because $\cu$ is a standard
iteration tree on the iterable model $\barw$.
The tree $\tcu$ also has a unique well founded branch $\tilde b$,
which must be 
the branch corresponding to $b$ since otherwise the preimage of
$\tilde b$ would be a second well founded branch through $\cu$.  Thus
we can define
$\pi_{\gth}\colon\cn_{\gth}\to\tcn_{\gth}$ to be the direct limit of
the maps $\pi_{\nu}\colon\cn_{\nu}\to\tcn_{\nu}$ for $\nu\in b$.
This
direct limit is defined since $\cu$ and $\tcu$ have the same
underlying tree and the maps $\pi_{\nu}$ commute with the respective
tree embeddings.

If $\phi$ is a limit ordinal, then we similarly have to pick a well
founded branch $b$ of $\ct\restrict\phi$.  
Since $\ct$ is not a standard iteration tree, we cannot use the
general theory as in the last paragraph.  Instead we use the following
lemma, the proof of which is deferred until after
the
construction of the trees $\ct$ and $\cu$.
\begin{lem}
There is a function $\mathfrak b$ such that if $\ct$ is defined at limit
ordinals $\phi$ by setting $[0,\phi]^{\ct}=\mathfrak b(\ct\restrict\phi)$
then every model $\cm_{\nu}$ of $\ct$ is well founded.
\label{l.strategy}
\end{lem}
\end{case}

\smallskip
This concludes case~1.  In the remaining cases both $\phi$ and
$\gth$ are successor 
ordinals, say $\phi=\gga+1$ and $\gth=\gga'+1$.  Let $\ga$ be the
largest ordinal such that $\cm_{\gga}$ and $\cn_{\gga'}$ agree up to
$\ga$, and let $\tilde\ga=\sup\pi_{\gga'}\image\ga$ and
$\cR=\ult(\cm_{\gga},\pi_{\gga'},\tilde\ga)$.

\begin{case}{2}{$\cm_{\gga}$ is a weasel and the phalanx
$\ppstree{\tilde\ga}{\tcu}\cR$ is not iterable}
This is the case in which the definition of $\ct$ differs from a
normal comparison.  We will define $\cm_{\phi}$ to be a premouse of
cardinality less than $\gd$, thus ensuring by brute force that if
$\phi$ is on the main branch of $\ct$ then the iteration along this
main branch will yield a set of indiscernibles of cardinality $\gd$.
The
following lemma will be proved later along with lemma~\ref{l.iterable2}:
\begin{lem}
If there is a ill behaved tree on $\ppstree{\tilde\ga}{\tcu}\cR$, then there is
an elementary substructure   $Q$ of
$\cm_{\phi}$, with $\ga\subset Q$ and $\card {Q}=\card\ga$, such that
there is an ill behaved tree on
$\ppstree{\tilde\ga}{\tcu}{\ult(Q,\pi_{\gth},\tilde\ga)}$. 
\label{l.iterable1}
\end{lem}
The case hypothesis asserts that the hypothesis of the lemma is true,
and we  
define $\cm_{\phi}$ to be the transitive collapse of
the elementary substructure
$Q$ of $\cm_{\gga}$ given by the lemma.  We extend the tree ordering
by letting $\nu\tlt{\ct}\phi$ if and only if $\nu<\phi$ and
$\nu\tlt{\ct}\gga$.
We will
regard this as a drop and hence leave the embedding $i_{\nu,\phi}$
undefined for $\nu\tlt{\ct}\phi$.  The ordinal $\rho_{\phi}$ 
associated with this stage of the tree (which
would be $\len(E_{\phi})$ in the normal successor case)
is defined to be the larger of 
$\sup_{\nu<\phi}\rho_{\nu}$ and the least ordinal $\gb$ such that
there is an ill behaved tree on the phalanx
$\ppstree{\bar\gb}{\tcu}{\ult(\cm_{\phi},\pi_{\gth},\bar\gb)}$, where
$\bar\gb=\sup\pi_{\gth}\image\gb$.

We will call the node $\phi$ of $\ct$ a {\em special successor node},
and we will call $\{\gga,\phi\}=\{\gga,\gga+1\}$ a {\em special pair}. 
Notice that $\gg$ and $\gg+1$ have the same set of predecessors in $\ct$,
even if $\gg$ is a limit ordinal.  This is impossible in a standard
iteration tree.
\end{case}

\begin{case}
{3}{Neither  case~1 nor case~2 holds}
This case is completely standard.  
Let $\ce_{\gga}$ and $\cf_{\gga'}$ be the extender sequences in 
 $\cm_{\gga}$ and $\cn_{\gga'}$ respectively, and
let $\beta$ be the least ordinal such that
$\gb\in\domain(\ce_{\phi})\cup\domain(\cf_{\gth})$, and 
$(\ce_{\phi})_{\beta}\not=(\cf_{\gth})_{\beta}$
if $\gb\in\domain(\ce_{\phi})\cap\domain(\cf_{\gth})$.
If $\gb$ is in the domain of $\ce_{\gga}$ then extend $\ct$ as follows: let
$\nu$ be the least  ordinal such that $\rho_{\nu}>\crit(E)$, let
 $\cm_{\gga}^*$ be the largest
initial segment of 
$\cm_{\nu}$ such that $E$ is an extender on $\cm_{\nu}$, and set
$\cm_{\gth}=\ult(\cm_{\gga}^*,E)$.
Similarly if 
$\cb$ is in the domain of 
$\cf_{\gga'}$ then extend $\cu\restrict\gth$
by setting  $\cn_{\gth}=\ult(\cn_{\gth}^*,F)$ and use the
shift lemma to define $\tcn_{\gth}=\ult(\tcn^*_{\gth},\pi_{\gth}(F))$
in $\tcu$, where $\tcn^*_{\gth}=\tcn_{\nu}$ if
$\cn^*_{\gth}=\cn_{\nu}$, and $\tcn^*_{\gth}=\pi_{\nu}(\cn^*_{\gth})$
if $\cn^*_{\gth}$ is a proper initial segment of $\cn_{\nu}$.
\end{case}
\subsection*{Proof of lemmas~\ref{l.strategy} and~\ref{l.iterable1}}
This completes the construction of the trees $\ct$, $\cu$ and $\tcu$
except for the proof of lemmas~\ref{l.strategy} and~\ref{l.iterable1}.
\begin{proof}[Proof of lemma~\ref{l.strategy}]
We  define an auxiliary iteration tree $\bct$, together with an
embedding from $\ct$ to $\bct$.  Since $\bct$ is a standard iteration
tree on $W$, all of its models are well founded. 
We show that at each limit ordinal $\phi<\len(\ct)$ there
is a unique branch  of $\ct$ which maps to the well founded branch
of $\bct$, and we will take $\mathfrak b(\ct\restrict\phi)$ to be this
unique branch. 

The embedding consists of a  map
$\gs\colon\len(\ct)\to\len(\bct)$, together with  embeddings
$t_{\nu}\colon\cm_{\nu}\to\bcm_{\gs(\nu)}$ for each ordinal $\nu<\len(\ct)$.   
The map $\gs$ is one to one and order preserving, with the
exception that both 
members of a special pair of $\ct$ map to the same node  of $\bct$.  
The construction is by recursion on the nodes $\phi$ of $\ct$.
Suppose that the construction has been carried out up to $\phi$, so we
have already constructed
$\bct\restrict\phi$ together with the map $\gs\restrict\phi\colon\phi\to\gth$
and the maps 
$t_{\ga}$ for $\ga<\phi$, and assume that the following conditions
are satisfied for all nodes $\ga,\gb<\phi$:

\begin{enumerate}
\item $\gs(\ga)\le_{\bct}\gs(\gb)$ if and only if either $\ga\le_{\ct}\gb$
or else there is $\ga'\le_{\ct}\gb$ such that
 $\{\ga,\ga'\}$ is a special pair.
\item If $\ga<\gb$ then
$t_{\ga}\restrict\rho_{\ga}=t_{\gb}\restrict\rho_{\ga}$. 
\item If there is a normal drop in $[0,\ga)_{\ct}$, so that
$\cm_{\ga}$ is an 
iterate of a mouse, then $t_{\ga}$ is a weak
$\deg(\ga)$-embedding in the sense of
\cite{mitchell-steel.inner-model}, and otherwise (if 
$\cm_{\ga}$ is a weasel or else the only drop has been at a special
successor) $t_{\ga}$ is an elementary embedding.
\end{enumerate}

\begin{case}
{1}{$\phi=\gga+1$ is a standard successor node of $\ct$}
This case is standard.  Set $\gs(\phi)=\gs(\gga)+1$, set
$\bare_{\gs(\gga)}=t_{\gga}(E_{\gga})$ and set
$\bcm_{\gs(\phi)}=\ult\left(\bcm^*_{\gs(\ga)},\bare_{\gs(\ga)}\right)$ where
$\bcm^*_{\gs(\ga)}$ is defined as usual.
Standard arguments show that the induction hypotheses are  true at $\phi$.
\end{case}

\begin{case}
2 {$\phi=\gga+1$ is a special successor node of $\ct$}
In this case we
set $\gs(\phi)=\gs(\gga)$.  By the construction of $\ct$ there is an
elementary embedding $s_{\gg}\colon\cm_{\phi}\to\cm_{\gga}$, so we can set
$t_{\phi}=t_{\gga}\cmp s_{\gg}$.  
Again, it is straightforward to verify that the induction hypotheses
are true at $\phi$.
\end{case}

\begin{case}3{$\phi$ is a limit ordinal}
Since clause~(1) holds for $\ga,\gb<\phi$ we have
$\phi=\sup\gs\image\phi$, so we can set 
$\gs(\phi)=\phi$.  Since $\bct\restrict\phi$ is a standard iteration
tree there is 
a unique well founded branch $\bar b$ through $\bct\restrict\phi$.
Extend $\bct$ as usual by letting
$\bcm_{\nu}=\bcm_{\bar b}$, the direct limit of $\bct$ along the
branch $\bar b$.  Now set 
$$b=\union\set{[0,\nu)_{\ct}:\gs(\nu)\in\bar b}
   = \set{\nu:\exists\nu'\in\gs^{-1}(\bar b)\;(\nu<_{\ct}\nu')}.$$ 
First notice that $b\subset\gs^{-1}(\bar b)$, since if $\nu\in b$ then
there is $\nu'\in\gs^{-1}(\bar b)$ such that $\nu<_{\ct}\nu'$, so
$\gs(\nu)<_{\bct}\gs(\nu')\in \bar b$. To see that $b$ is linearly
ordered, suppose $\nu_0,\nu_1\in b$ with $\nu_0<\nu_1$.  If
$\nu_0\not<_{\ct}\nu_1$ then clause~(1) implies that there is $\nu_0'$
such that $\set{\nu_0,\nu_0'}$ is a special pair and
$\nu_0'\le_{\ct}\nu_1$.  Since $\nu_0\in b$ there is $\nu_1'$ such that
$\gs(\nu_1')\in \bar b$ and $\nu_0<_{\ct}\nu_1'$, but this is
impossible because in this case $\gs(\nu_1)$ and $\gs(\nu'_1)$ would
be incomparable in~$\bct$.

Hence $b$ is a cofinal branch of $\ct\restrict\phi$, so we can define
$\cm_{\phi}=\cm_b$, the limit of $\ct$ along the branch $b$.  Now define
$t_{\phi}$ to be the limit of the embeddings $t_{\nu}$ for $\nu\in b$.
Again, it is straightforward to verify the induction hypotheses.
\end{case}
This completes all three cases, and hence the proof of
lemma~\ref{l.strategy}. 
\end{proof}

We are now ready to to finish part one of the proof by proving
lemma~\ref{l.iterable1}.
The proof is essentially the same as that of the following lemma, 
which is needed for part two of
the proof, so we will combine the two proofs.
This technique of applying Martin-Steel iterability theorem 
\cite{martin-steel.iteration-trees} is due to Woodin. 

\begin{lem}
Suppose $\cm_{\phi}$ and $\cn_{\gth}$ are the last models
of $\ct$ and $\cu$, so that 
$\cm_{\phi}$ and $\cn_{\gth}$ are matched
up to $\gd$, and suppose that  
$\cm_{\phi}$ is a weasel.  Then 
there are no ill behaved trees on
the phalanx
$\ppstree{\tilde\gd}{\tcu}{\ult(\cm_{\phi},\pi_{\gth},\tilde\gd)}$,
where
$\tilde\gd=\sup\pi_{\gth}\image\gd$.
\label{l.iterable2}
\end{lem}
\begin{proof}
[Proof of lemmas~\ref{l.iterable1} and~\ref{l.iterable2}]
For lemma~\ref{l.iterable1}
we are given models $\cm_{\phi}$ and $\cn_{\gth}$ from the trees $\ct$
and $\cu$ and an ordinal $\ga<\gd$ such that $\cm_{\phi}$ and
$\cn_{\gth}$ agree up to $\ga$.  
In the case of lemma~\ref{l.iterable2} we take $\cm_{\phi}$ and
$\cn_{\gth}$ to be  the final models of the trees $\ct$ and $\cu$, and
we take $\ga=\gd$.

 The map $\pi$ was used to copy the tree $\cu\restrict\gth$ on
$L_{\eta}[\ce']$ to a tree $\tcu$ on $L[\ce]$, with copy 
map $\pi_{\nu}\colon \cn_{\nu}\to \wt \cn_{\nu}$ for each
model $\cn_{\nu}$ of $\cu$. 
 Set
$\tilde\ga=\sup\pi\image\ga$ and
$\cR=\ult(\cm_{\phi},\pi_{\gth},\tilde\ga)$, and assume 
that there is a ill behaved tree on
$\ppstree{\tcu}{\tilde\ga}{\cR}$. 
In case of lemma~\ref{l.iterable1} the tree is part of the
hypothesis, while for lemma~\ref{l.iterable2} we assume the
existence of such an ill behaved tree towards a proof by
contradiction. 

Let $Z_1\prec\ps(H_\gl )$ be an elementary
substructure such that $\{\ca,C,W,X\}\subset Z_1$, and such that
\begin{itemize}
\item For lemma~\ref{l.iterable1}: $\ga+1\subset Z_1$ and
$\card{Z_1}=\card\ga$.  We set $\eta=\ga$.
\item For  lemma~\ref{l.iterable2}: $Z_1\cap\gd=\sup(Z_1\cap\gd)=\card
{Z_1}$.   We set $\eta=Z_1\cap\gd$. 
\end{itemize}
Notice that everything which has been mentioned is definable from
members of $Z_1$, and hence is in $Z_1$.
We set $\tilde\eta=\sup\pi\image\eta$.

Now let $\psi_1\colon M_1\cong Z_1$ be the transitive collapse,
and let $Q=\psi_1^{-1}(\cm_\phi)\cong \cm_{\phi}\cap Z_1$.
The major part of the proof of these lemmas is the proof of the
following claim:
\begin{claim}
There is an ill founded tree $S_1$ on
$\pstree{\tilde\eta}{\tcu}{\ult(Q,\pi_{\gth},\tilde\eta)}$. 
\label{c:itA}
\end{claim}
In the case of lemma~\ref{l.iterable1}, this claim is exactly what
is required.  
Before proving claim~\ref{c:itA} we show that
lemma~\ref{l.iterable2} also follows from this claim.

\begin{proof}[Proof of \ref{l.iterable2} from \ref{c:itA}]
Note that $\eta\in [0,\phi)_{\ct}$ since $[0,\phi)_{\ct}$ is
closed and unbounded in~$\gd$.  We claim that $Z_1\cap
\cm_{\phi}\subset\range(i_{\eta,\phi})$, so that 
the embedding $\psi_1\colon Q\cong Z_1\cap \cm_{\phi}\prec \cm$ can be
factored  
\[
\psi_1\colon Q
\xrightarrow{k} \cm_{\eta}
\xrightarrow{i_{\eta,\phi}^{\ct}}\cm_{\phi}
\]
where $k=(i^{\ct}_{\eta,\phi})^{-1}\cmp\psi_1$.  

To see that  $Z_1\cap M_{\phi}\subset\range(i_{\eta,\phi})$, first
  notice that  every member of 
$\cm_{\phi}$ is in $\range(i^{\ct}_{\nu,\phi})$ for some 
$\nu\in [0,\phi)_{\ct}$.  Hence any member $x$ of $M_\phi\cap Z_1$ is in 
$\range(i^{\ct}_{\nu,\phi})$ for some    
$
\nu\in [0,\phi)_{\ct}\cap   
Z_1=[0,\phi)_{\ct}\cap\eta=[0,\eta)_{\ct},  
$
since $\eta\in [0,\phi)_{\ct}$,  
and hence $x\in \range(i^{\ct}_{\eta,\phi})$.  

The embedding $k\colon Q\to \cm_{\eta}$ induces a map
\[
\tilde k\colon
\ult(Q,\pi_{\gth},\tilde\eta)\to\ult(\cm_{\eta},\pi_{\gth},\tilde\eta),
\]
which in turn can be used  
to copy $\cs_1$ to a tree $\cs_2$ on 
$\pstree{\tilde\eta}{\tcu}{\ult(\cm_{\eta},\pi_{\gth},\tilde\eta)}$ 
which is necessarily also ill behaved.  Now 
$\pi_{\gth}\restrict\eta=\pi_{\gth'}\restrict\eta$, where $\cn_{\gth'}$
is the stage which the tree $\cu$ had reached at the time $\cm_{\eta}$
was being considered in tree $\ct$.   
Thus
$\pstree{\tilde\eta}{\tcu}{\ult(\cm_{\eta},\pi_{\gth},\tilde\eta)}$ 
is the tree $\pstree{\tilde\ga}{\tcu}\cR$ of case~(2) of the
construction of the trees, so the existence of the ill behaved tree $\cs_2$ 
would have
caused the tree $\ct$ to drop at stage $\eta$.  Furthermore
$\rho_{\eta}=\eta=\union_{\nu<\eta}\rho_{\nu}$, so the final model
$\cm_{\phi}$ of
the tree must be above the second member $\eta+1$ of the special pair,
rather than above $\eta$.
This contradicts the assumption that $\cm_{\phi}$ is a weasel and
hence completes the proof of lemma~\ref{l.iterable2} 
from claim~\ref{c:itA}.
\end{proof}

\begin{proof}
[Proof of claim~\ref{c:itA}]
To find $\cs_1$, let $Z_0\prec Z_1$ with $\set{A,C,W,X}\subset Z_0$
and $\card {Z_0}=\gw$.  Then $Z_0$ satisfies that there is an
ill behaved tree $\cs$ on 
$\pstree{\tilde\ga}{\tcu}{\cR}$.  
Let $M_0\cong Z_0$ be transitive, with maps
\[
\psi\colon M_0\xrightarrow{\psi_0} M_1 
             \xrightarrow{\psi_1}\ps(H_\gl ).
\]
Set $\cs_0=\psi^{-1}(\cs)$.  Then $M_0$ satisfies that 
$\cs_0$ is ill behaved, but we need to show that $\cs_0$ has no
branches even in $V$:
\begin{claim}
The tree $\cs_0$ is ill behaved in $V$.\label{c:itB}
\end{claim}
\begin{proof}
Suppose to the contrary that there 
is a well founded branch through $S_0$ in $V$.

Let $\cf$ be the limit of the extender sequences of the models in
$\cs$; that is, if $\xi<\len(\cf)$ then $\cf\restrict\xi$  is an
initial segment of the $\nu$th model $\cm^{\cs}_\nu$ for
every sufficiently large $\nu<\len(\cs)=\eta$; but $\cf$ itself is
not an 
initial segment of the extender sequence of any of these models.
By
our assumption there is no class model with a Woodin cardinal, and
hence there is $\gga$ such that $L_{\gga}[\cf]$ satisfies that
$\len(\cf)$ is not a Woodin cardinal.  By elementarity $\gg\in M_0$.

Set $\gga'=\psi^{-1}(\gga)$, and let $G$ be
$M_0$-generic for the Levy collapse $\col(\gw,\gga')$ of $\gga'$ onto
$\gw$.    In $M_0[G]$ form 
the tree $A$ of attempts to find a branch through $\cs_0$ such that the
limit along the branch either is well founded or has a well founded
part of length 
at least $\gga$.  Because $\cs_0$ has a well founded branch 
in $V$, the tree $A$ has an infinite branch in $V$, and  by the
absoluteness of well order it follows that $A$ has an infinite
branch in $M_0[G]$ as well.  
Thus $\cs_0$ has a well founded branch~$b$ in~$M_0$.
By the 
Martin-Steel iterability  theorem
\cite{martin-steel.iteration-trees} there can be at most
one branch $b\in M_0[G]$ through $\cs_0$ which is well founded up to
$\gga$. 
The uniqueness of $b$, together with  the homogeneity
of  the Levy collapse  $G$, implies that $b\in M_0$.
This contradicts the fact that $M_0$ satisfies that $\cs_0$ has no
well founded branch, and hence completes the proof of
claim~\ref{c:itB}. 
\end{proof}
In order to complete the proof of claim~\ref{c:itA}, and hence of
lemmas~\ref{l.iterable1} and~\ref{l.iterable2}, we  copy
$\cs_0$ to a tree on 
$\ppstree{\tilde\eta}{\tcu}{\ult(Q,\pi_{\gth},\tilde\eta)}$ as follows:

The tree $\cs_0$ is on the phalanx
\begin{equation}
\psi^{-1}\left(
\pstree{\tilde\ga}{\tcu}{\ult(\cm_{\phi},\pi_{\gth},\tilde\ga_1)}\right).
\label{ph1}
\end{equation}
If $P_\nu$ is the $\nu$th member of the phalanx~\eqref{ph1} then 
$\psi_{0}\restrict P_\nu\colon P_{\nu}\to\psi_{0}(P_{\nu})$,
which is the $\psi_{0}(\nu)$th member of the phalanx
\begin{equation}
\psi_{1}^{-1}\left(\pstree{\tilde\ga}{\tcu}
       {\ult(\cm_{\phi},\pi_{\gth},\tilde\ga)}\right)
=
\pstree{\tilde\ga_1}{\psi_1^{-1}(\tcu)}
       {\ult(Q,\psi_1^{-1}(\pi_{\gth}),\ga_1)}
\label{ph2}
\end{equation}
where $\tilde\ga_1=\psi^{-1}(\tilde\ga)$.
The phalanx~\eqref{ph2} can embedded into
$\pstree{\tilde\eta}{\tcu}{\ult(Q,\pi_{\gth},\tilde\eta)}$ by the map
$\psi_1$, and by copying $\cs_0$ along this embedding we 
obtain the required  ill behaved tree on  
$\pstree{\tilde\eta}{\tcu}{\ult(Q,\pi_{\gth},\tilde\eta)}$.
\end{proof}
This completes the proof of claim~\ref{c:itA} and hence of
lemmas~\ref{l.iterable1} and~\ref{l.iterable2}.
\end{proof}

\subsection*{The Indiscernibles Generated by $\ct$}
This concludes the construction of the trees $\ct$ and $\cu$, which is
the first half 
of the proof of theorem~\ref{t:j->erdos}, and we are now
ready to use $\ct$ and $\cu$
to show that $L[\ce]$ has a set of indiscernibles for the
structure $\ca$.
First, we will show that the tree $\ct$ has length~$\gd$ and hence
generates a cofinal set of indiscernibles for its last model.  Using
this, we will then show that the tree $\cu$ does not drop along its
main branch.  This implies  that the set $I$ of indiscernibles
from $\ct$ 
is a set of 
indiscernibles for every set in $\barw$, and in particular for
$\pi^{-1}(\ca)$.  Finally we will show that the filter on $L[\ce]$
generated by $\pi\image I$ is a member of $L[\ce]$, and use this
filter inside $L[\ce]$ to construct the required set $D\in L[\ce]$ of
indiscernibles for $\ca$.  The 
next lemma, which shows that $\ct$ generates a set of indiscernibles
of size $\gd$, is the main lemma of this half of the proof.

\begin{lem}
The tree $\ct$ has length $\gd+1$, and there are ordinals $\nu<\gd$ and
$\rho<\gd$ such that $i_{\nu,\gd}$ is  defined and $i_{\nu,\gd}(\rho)\ge\gd$.
\label{l.delta-indisc}
\end{lem}
\begin{proof}
Let $\phi+1$ be the length of the \qi{} tree $\ct$ on $L[\ce]$, and let 
 $\gth+1$ be the length of the iteration tree $\cu$  on
 $\barw=\pi^{-1}(L[\ce])$.  
Since $\gd$ is a cardinal, the proof of the comparison lemma implies
that $\phi,\gth\le\gd$.  

\begin{claim}
The final models $\cm_{\phi}$ and $\cn_{\gth}$ of the trees $\ct$ and
$\cu$ have size at least
$\gd$, and hence  agree up to $\gd$.
\end{claim}
\begin{proof}
Suppose to the contrary that one of these models has size less than
$\gd$.  By the construction of the trees $\ct$ and $\cu$ it follows
that that model is an initial segment of
the final model of the other tree.  In addition, that tree must drop
along its main 
branch, since the each of the roots $\barw$ and $W$ of the two trees
 have
 cardinality at least $\gd$. 
Lemma~\ref{l.compare}(1) asserts that this is impossible in any
comparison using ordinary iteration trees.  We will verify that this
is still impossible for the modified comparison using the \qi{} tree
$\ct$.   

The reason why lemma~\ref{l.compare}(1) is true is that if there is an
(ordinary) drop at some node $\nu+1$ in the main branch of one of the
trees, say $\ct$, then 
there is a subset $x$ of the projectum $\rho$ of $\cm_{\nu+1}$
which is 
definable in $\cm^{\ct}_{\nu+1}$ and is not a member 
of the corresponding model of the other tree.  If we take $\nu+1$ to
be last place along the main branch where such a drop occurs, then
$i_{\nu+1,\gth}^{\ct}$ exists, and $i_{\nu+1,\gth}^{\ct}\restrict\rho$
is the identity.  Then $x$ is definable in $\cm_{\gth}$ and not a
member of $\cn_{\phi}$, so $\cm_{\gth}$ cannot be a proper initial
segment of $\cn_{\phi}$.

This argument is unaffected by the use of the \qi{} tree $\ct$ instead
of a normal iteration tree.  It follows that $\ct$ must have a special
drop $\nu+1$, but no normal drops, on its main branch, and
$\cm_{\gth}$ must be an initial segment of $\cn_{\phi}$.   From
the definition of a special drop it follows that the phalanx
$\cu\oplus_{\tilde\rho}\ult(\cm_{\nu+1},\pi_{\nu'},\tilde\rho)$ is not
iterable, where
$\nu'$ is the stage in $\cu$ corresponding to $\nu+1$ in $\ct$, and
$\tilde\rho=\rho^{\ct}_{\nu+1}$.   We will use this lack of
iterability just like the set $x$ in the case of a normal drop: it
implies that the phalanx
$\cu\oplus_{\tilde\gd}\ult(\cm_{\gth},\pi_{\phi},\tilde\gd)$ is not
iterable, but 
if $\cm_{\gth}$ is an initial segment of $\cn_{\phi}$ then
$\ult(\cm_{\gth},\pi_\phi,\tilde\gd)$ can be embedded into
$\ult(\cn_{\gth},\phi_{\phi},\tilde\gd)$, so that the phalanx 
$\cu\oplus_{\tilde\gd}\ult(\cn_{\phi},\pi_{\phi},\tilde\gd)$ is not
iterable.  This is absurd, since the later phalanx is actually the
standard iteration tree $\tcu$ on the iterable model $W$.

This contradiction completes the proof of the claim.
\end{proof}

Thus the models $\cm_{\phi}$ and $\cn_{\gth}$ agree up to $\gd$.  
If there is a drop of any kind in the main branch of
$\ct$  then $\card{\cm_{\nu}}<\gd$ for every sufficiently large $\nu$
in the main branch of $\ct$.
Since $\card{\cm_{\phi}}=\gd$  while
all of the extenders in $\ct$ have length less than $\gd$,
it follows that $\phi=\gd$ as
required.  

\newcommand{\othermap}{t}
\newcommand{\Rtree}{\pstree{\tilde\gd}{\tcu}{\cR}}
\newcommand{\rtree}{\Rtree\pcat\cw}
\newcommand\utree{\tcu\pcat\cv}

Thus we can assume for the remainder of the proof of this lemma that
there are no 
drops in the main branch of $\ct$, and hence $\cm_\phi$ is a weasel.
If we set 
$\tilde\gd=\union\pi_{\gth}\image\gd$ and
$\cR=\ult(\cm_{\gd},\pi_{\gth},\tilde\gd)$ then lemma~\ref{l.iterable2}
asserts that there are no ill behaved trees on the phalanx
 $\Rtree$, so we can compare the
models $\cR$ and $\tcn_{\gth}$, using an iteration tree $\rtree$ on the phalanx
$\Rtree$ and an iteration tree $\utree$ on the
phalanx $\tcu$.  The comparison takes place as if $\tcu$ were simply a
phalanx; however we will later make use of the iteration tree
structure on $\tcu$, regarding $\rtree$ as
a (nonstandard) iteration tree with two roots $L[\ce]$ and $\cR$.

\begin{claim}
\enumerate
\item There are no drops on the main branch of either of the trees
$\rtree$ or $\utree$.
\item The trees $\rtree$ and $\utree$ have the same last model.
\item
The final model of $\rtree$ is above $\cR$, rather than above $\tcn_0$.
\label{c.di}
\end{claim}
\begin{proof}
The proof of this claim uses the techniques of
lemma~\ref{l.phi-minimal}(5).  The proofs of clauses~(1)
and~(2) are the same as the proof of the corresponding facts in
lemma~\ref{l.phi-minimal}(5), using the fact that  
there is no drop in $\ct$ and hence there
is an elementary embedding from $L[\ce]$ into $\cR$.  It follows
that $\cR$
is universal if hypothesis~(1) of the main theorem is true, and $\cR$ is
$\varphi$-minimal if hypothesis~(2) of the main theorem is true.
If there is no drop in the main
branch of $\rtree$ then the same will be true of the final model $P$
of $\cw$, 
since there is an elementary embedding from either $L[\ce]$ or $\cR$
into $P$.

The proof of clause~(3) is similar to the proof in
lemma~\ref{l.phi-minimal}(5) of the fact that $Q$ is not above
$L[\ce]$. 
Suppose to the contrary that $P$ lies above
the root $\tcn_0=L[\ce]$ of $\tcu$
in the tree $\rtree$. Let $b$ be the main branch of $\rtree$, 
let $c$ be the main branch of $\utree$
and let $\nu\in\domain\cu$ be the largest member of $b\cap c$.   
Finally, let $E$ be the first extender used along $b$ after $\nu$, and
let $E'$ 
be the first extender used along $c$ after $\nu$.  
Then we can  write 
the embedding $i_b$ along the main branch $b$ of $\rtree$ in the form
\[
\begin{CD}
  i_b\colon
  L[\ce]@>{i^{\cu}_{0,\nu}}>>\tcn_{\nu}@>{i^{E}}>>\ult(\tcn_{\nu},E)@>k>>P,  
\end{CD}
\] 
and write the embedding $i_c$ along the main branch of $\utree$ in the form
\[\begin{CD}
i_c\colon
L[\ce]@>i^{\cu}_{0,\nu}>>\tcn_{\nu}@>i^{E'}>>\ult(\tcn_{\nu},E')@>k'>>P.
\end{CD}
\]

By lemma~\ref{l.phi-minimal}(4) or lemma~\ref{l.steel-core-model}(4)
we get $i_b\restrict K_\gk=i_{c}\restrict K_{\gk}$.
Let $\eta$ be $\min(\crit(E), \crit(E'))$.  Then since $\eta<\gk$ it follows
that $k\cmp i^{E}\restrict\ps(\eta)=k'\cmp i^{E'}\restrict\ps(\eta)$ and hence
one of $E$ and $E'$ is an initial segment of the other.  By the proof
of lemma~\ref{l.compare}(\ref{l.comp3}) this can never happen if either of the
extenders $E$ or
$E'$ come from $\cw$ or $\cv$, so both of $E$ and $E'$ must come from 
$\cu$. This implies that $E=E'$, contradicting the
choice of $\nu$.
This completes the proof of the claim.
\end{proof}

The following diagram illustrates the present situation, with the
straight arrows indicating embeddings which are known to exist:
\begin{lrbox}
{\diagrambox}
\xymatrix@C=.5pc{
       &&{P}&{P}\\
{\cm_{\phi}} \ar[rr]_{\tilde\pi_{\gth}}
       &&{\cR} \ar[u]_{i^{\cw}}\\
{L[\ce]}  \ar[u]_{i^{\ct}}
       &&{L[\ce]} \ar@{~>}[u]_{\cu}
           &{L[\ce]} \ar[uu]_{s=i^{\tcu\oplus\cv}}
}
\end{lrbox}
\begin{equation*}
\usediagram
\end{equation*}
The left hand map $i^{\ct}$ is defined since there are no
drops on the main branch of $\ct$.  Thus there are two maps from
$L[\ce]$ to $P$, namely $s=i^{\tcu\oplus\cv}$ and
$\othermap=i^{\cw}\cmp\tilde\pi_{\gth}\cmp i^{\ct}$.  Then
$\othermap\restrict L_{\gl}[\ce]= s\restrict L_{\gl}[\ce]$ by
lemma~\ref{l.phi-minimal} or~\ref{l.steel-core-model}, and in
particular $s$ and $t$ have the same critical point $\rho$, with
$\othermap(\rho)=s(\rho)$ and
$\othermap\restrict\ps(\rho)=s\restrict\ps(\rho)$.

We claim that $s(\rho)\ge\tilde\gd$.  Suppose,  
 to the contrary,	that  
$s(\rho)<\tilde\gd$ and let $E$ be the first extender used in the
main branch of $\utree$.  Then $\rho=\crit(E)$, and  
$\tilde\gd >i^{E}(\rho)>\len(E)$ so $\len(E)<\tilde\gd$.  It
follows 
that $E$ comes from $\tcu$, so that $E=\pi_{\nu}(\bare)$ for some
$\nu\le\gth$, where $\bare=E^{\cu}_{\nu}$.  Now let 
$\xi=\crit(\pi)=\crit(\pi_{\nu})=\crit(\pi_{\gth})$.  Then 
$\rho=\crit(t)\le\xi$, and it follows that $\rho<\xi$ since
$\rho=\crit(E)=\pi(\crit(\bare))$ and $\xi\notin\range(\pi)$.
Since $\rho=\crit(t)<\crit(i^{\cw}\cmp\pi_{\gth})$, it follows that  
$\rho=\crit(i^{\ct})$ and hence the main branch of $\ct$ begins 
with 
an extender $F$ such that $\crit(F)=\rho$. Since
$\othermap\restrict\ps(\rho)=s\restrict\ps(\rho)$  
it
follows that one of $F$ and $\bare$ is an initial segment of the
other, contradicting the construction of $\ct$ and~$\cu$.

Thus $ s(\rho)\ge\tilde\gd$, and since $t(\rho)=s(\rho)$ we also have
$t(\rho)=i^{\cw}\cmp\pi_{\gth}\cmp i^{\ct}(\rho)\ge\tilde\gd$.
Since $i^{\cw}\restrict\tilde\gd$ is the 
identity it follows that $\pi_{\gth}\cmp
i^{\ct}(\rho)\ge\tilde\gd$, and 
hence $i^{\ct}(\rho)\ge\gd$.  But each of the extenders in the tree $\ct$
has length less than $\gd$, so the length of the main branch of $\ct$,
and hence of $\ct$ itself, cannot be smaller than~ $\gd$.
This completes the proof of lemma~\ref{l.delta-indisc}.
\end{proof}
Now we can look at the indiscernibles generated by $\ct$.

\begin{cor}
There is no drop in the main branch of $\cu$.  Furthermore
$j_{0,\gth}\image\gd\subset\gd$ and there is a closed and
unbounded  subset $I$
of the main branch $[0,\gd]_{\ct}$ of $\ct$ satisfying the following
four conditions:
\label{c:indisc}
\begin{enumerate}
\item If $\nu\in I$ then $\nu=\crit(i_{\nu,\phi})$.
\item If $\nu,\nu'\in I$ with $\nu<\nu'$ then
$\nu'=i_{\nu,\nu'}(\nu)$.
\item If $\nu\in I$ then $j_{0,\gth}\image\nu\subset\nu$.
\item Every member of $\nu$ is regular in $\barw$, and is a limit
member of $\pi^{-1}(C)$.
\end{enumerate}
\end{cor}
\begin{proof}
If there were a drop on the main branch of $\cu$ then there would be a
model $\cn_\nu$ on the main branch of $\cu$ with
$\card{\cn_\nu}<\gd$.  If there is also a drop on the main branch of
$\ct$ then the proof of the comparison lemma shows that both trees
have length less than $\gd$; while if there is no drop on the main
branch of $\ct$ then lemma~\ref{l.delta-indisc}
implies that $\crit(E_{\nu})\le i_{0,\nu}(\rho)$ for every extender
$E_\nu$ used on the main branch of  
$\ct$, and again it follows by the proof of the comparison lemma that both
trees have length less than $\gd$.

Since either case contradicts lemma~\ref{l.delta-indisc}, it
follows that there is no drop on the main branch of $\cu$. 
The rest of the proof of corollary~\ref{c:indisc} is straightforward.
\end{proof}
The indiscernibles $I$ are not true indiscernibles, since they come from
ultrapowers by different ultrafilters.  
If $M$ is any  model of the form $L[\ce]$ for some good sequence $\ce$
then define
\begin{equation*}
M\models \nu\in_0 x\iff
  \begin{cases}
    \nu\in x &\text{if $o^{M}(\nu)=0$}\\
    x\cap\nu\in U_\nu &\text{if $o^M(\nu)>0$}
  \end{cases}
\end{equation*}
where $U_\nu$ is the unique order~$0$ measure on $\nu$ in $M$.

Define  $\baru$  to be the set of subsets $x$ of
$\gd$ such that $\barw\models\nu\in_0 x$ for all
$\nu\in I$.
We begin by showing that $\baru$ is a normal measure on $\barw$.  We
will then use this fact to show that the filter generated in the same
way by
$\pi\image I$ is an ultrafilter on the model $\ca$ for which we need
to find indiscernibles.

\begin{defn}
If $F$ is a filter on $\ps(\xi)$ for some ordinal $\xi$ then 
the \emph{$F$-closure} of a set $N$ is the smallest set
$Y\subset\bigcup_{n<\gw}\ps([\xi]^n)$ such that
$N\cap\bigcup_{n<\gw}\ps([\xi]^n)\subset Y$ and,
$\set{\vec\nu:\set{\ga:(\nu_0,\dots\nu_{n-1},\ga)\in x}\in F}\in
Y$ whenever $x\in Y$.
\end{defn}
\begin{lem}
Let  $ Y$ be the $\baru$-closure of 
$\barw$.
Then  $\baru$ is a normal ultrafilter on~ $Y$.
\end{lem}
\begin{proof}
First we will show that $\ps(\gd)\cap\cn_{\gth}\subset\cm_{\gd}$.
Suppose $x\in\cn_{\gth}$ and $x\subset\gd$.  Then
$x\cap\nu\in\cm_{\nu}$  
for every sufficiently large $\nu\in I$, and hence there is a
stationary subset $I'\subset I$ such that for all pairs $\nu<\nu'$ of
members of $I'$
we have $i_{\nu,\nu'}(x\cap\nu)=x\cap\nu'$. But
then $x=i_{\nu,\gd}(x\cap\nu)\in\cm_\gd$ where $\nu$ is any member  of $I'$.

Now let $U'$ be the filter defined like $\baru$, but using the 
order~0 measures $U'_{\nu}$ from 
$\cn_{\gth}$.
 That is, $x\in U'$ if and only if $\cn_{\gth}\models\nu\in_0x$ for
every sufficiently large $\nu\in I$. 
Then $U'$ is an ultrafilter on the
$U'$-closure of $\cm_{\gd}$ and hence on the $U'$-closure of
$\cn_{\gd}$. Now we claim that $\baru=\set{x:j_{0,\gth}(x)\in U'}$.  
Let $\nu$ be in $I$.  If $o^{\barw}(\nu)=\nu$ then
$j_{0,\gth}(\nu)=\nu$ by corollary~\ref{c:indisc}(3,4), and hence
\[\barw\models \nu \in_0 x\iff 
\nu\in x\iff\nu\in j_{0,\gth}(x)\iff \cm_{\gth}\models \nu\in_0
j_{0,\gth}(x). 
\] 
If $o^{\barw}(\nu)>0$ 
then the situation is slightly more complicated.   
We have $\barw\models\nu\in_0 x$ if and only if
$x\cap\nu\in\baru_\nu$, and  
if
$o^{\cn_{\gth}}(\nu)>0$ then 
\[
x\cap\nu\in\baru_{\nu}
\iff
j_{0,\gth}(x)\cap\nu\in U'_{\nu}
\iff
\cm_{\gth}\models\nu\in_0 j_{0,\gth}(x).
\]
On the other hand, if $o^{\cn_{\gth}}=0$ then $j_{\nu,\nu+1}\colon
\cm_{\nu}\to\ult(\cm_\nu,j_{0,\nu}(\baru_{\nu}))$ so  
\[
x\cap\nu\in\baru_{\nu}\iff
\nu\in j_{0,\gth}(x)
\iff
\cm_{\gth}\models\nu\in_0j_{0,\gth}(x).
\]
Since $U'$ is a normal ultrafilter on the $U'$-closure of $\cn_{\gth}$
it follows that $\baru$ is a normal ultrafilter on the $\baru$-closure
$Y$ of $\barw$.
\end{proof}

Now repeat the process in $W$, defining a filter $U$ on
$\tilde\gd=\sup\pi\image\gd$ by $x\in U$ if and only if
$W\models\nu\in_0x$ for every sufficiently large $\nu\in\pi\image I$.

The filter $U$ need not be an ultrafilter on $W$, but we will find a
premouse containing $\ca$ on which $U$ is an ultrafilter.
In order to do so let $\mse m\in X$ be the least mouse such that $\ca$
and $C$ are members of $\mse m$,  
let $h$ be the canonical
Skolem function of $\mse m$, and let $\mse m^*$ be the transitive
collapse of $h\image\tilde\gd$.  

\begin{lem}
\label{UisUF}
The filter
$U$ is a normal ultrafilter on  the $U$-closure of $m^*$.
\end{lem}
\begin{proof}
First we show that $U$ is an ultrafilter on $\mse m^*\cap\ps(\tilde\gd)$.
Let  $x$ be an arbitrary subset of  $\tilde\gd$ in $\mse m^*$.
Then there is an ordinal $\ga<\tilde\gd$ such that
$x=h(\ga)\cap\gd$. 
Now for sufficiently large $\nu,\nu'\in I$ we have  
\begin{equation}
\barw\models \forall\gb<\nu\,\bigl(\nu\in_{0}
\bar h(\gb)\iff{\nu'}\in_{0}\bar h(\gb)\bigr).
\label{e3}
\end{equation}
Pick $\nu_0\in I$ such that $\ga<\pi(\nu_0)$ and equation~\eqref{e3}
holds
 for all $\nu'>\nu\ge\nu_0$ in $I$.
Then 
\[
W\models \forall\gb<\pi(\nu)\,\bigl(\pi(\nu)\in_0
\bar h(\gb)\iff\pi(\nu')\in_0 h(\gb)\bigr),
\]
also holds for all $\nu'>\nu\ge\nu_0$.  In particular, 
$\pi(\nu)\in_0 h(\ga)\iff \pi(\nu')\in_0 h(\ga)$, holds for all such $\nu$
and $\nu'$, so that either $x\in U$ or
$\tilde\gd\setminus x\in U$.  Since $x$ was arbitrary it follows that
$U$ is an  
ultrafilter on $\mse m^*$. A straightforward
extension of this argument proves that $U$ is a normal ultrafilter
on the $U$-closure of $\ps(\tilde\gd)\cap\mse m^*$.
\end{proof}

\begin{lem}
\label{UinW}
The $U$-closure $Y$ of $\set{\smash{x\cap\tilde\gd:x\in \mse m^*}}$ is a member of
$L[\ce]$. 
\end{lem}
\begin{proof}
Let $\mse n$ be the least mouse such that there is a subset of
$\tilde\gd$, definable in $\mse n$, which is not measured by $U$; or
if there is no such mouse then let $\mse n=W$.  Then $U\cap\mse n$ is
a normal measure.  We will finish the proof of lemma~\ref{UinW} by
proving the following claim, which clearly implies lemma~\ref{UinW}:
\begin{claim} 1.~$Y\subset \mse n$, and 2.~$U\cap\mse n\in K$.
\label{c.1}
\end{claim}

For clause~(1), let $\mse m^*$ be as
above, and let $i^{U}_{n}\colon m^*\to\mse m_n=\ult_n(\mse
m^*,U)=\ult(\mse m^*,U^n)$ be the 
$n$-fold iterated ultrapower, which is defined by lemma~\ref{UisUF}.
Notice that every set in $\mse m_n$ is
measured by $U$, and that 
\begin{multline*}
\set{\nu_0:\set{(\nu_1,\dots,\nu_{n+1} ):(\nu_0,\dots,\nu_{n+1})\in
x\cap\gd}\in U^n}
\\
=\set{\nu_0<\tilde\gd:(\nu_0,\tilde\gd,i^{U}_1(\tilde\gd),\dots,i^U_{n-1}(\tilde\gd))
\in i^U_{n}(x)}\in \mse m_n
\end{multline*}
for any set $x\in X$.  Thus $Y\subset \union_{n}\mse m_n$, and
it will be enough to show that the subsets of $\tilde\gd$ in
$\mse m_n$ are an initial subset of those in $K$, so that $\mse m_n$
is an initial segment of $\mse n$.

For clause~(2), consider the model $\mse n_1=\ult(\mse n,U)$.  Since
$U\cap\mse n$ is definable from $\mse n_1$, it will be enough to show
that $\mse n_1\in W$.  

Both clauses follow from the following claim:

\begin{claim}
\label{c.2}
Suppose that $\mse p$ is a premouse which agrees with $\cw$ up to
$\tilde\gd$, that every member of $\mse p$
is 
definable from parameters in $\tilde\gd\cup p$ for some finite set $p$,
and that the phalanx $(\tilde\gd,(W,\mse p))$ is iterable.  Then
$\mse p\in W$ and $\ps^{\mse p}(\tilde\gd)$ is an initial segment of
$\ps^{W}(\tilde\gd)$.
\end{claim}

 Claim~\ref{c.1} will follow from claim~\ref{c.2}, provided we can
show that the phalanx $(\tilde\gd,(W,\mse n_1))$ and the phalanxes
$(\tilde\gd,(W,\mse m_n))$ for $n<\gw$ are all iterable. 
\begin{proof}[Proof of claim~\ref{c.2}]
Compare the given phalanx with $W$, using trees $\cs$ on
$(\tilde\gd,(W,\mse p))$ 
and $\ct$ on $W$.  
Standard arguments show that the last model of
$\cs$ must lie above $\mse p$, that there are no drops on the main
branch of $\cs$, and that the final model $P$ of $\cs$ is an initial
segment of the final model $Q$ of $\ct$.   The models $P$ and $\mse
p$ have the same subsets of $\tilde\gd$, so the subsets of $\tilde\gd$
in $\mse p$ are an initial segment of those in $Q$, and hence of those
in $W$.  Furthermore, $\mse p$ is isomorphic to the Skolem hull in $P$
of $\tilde\gd\cup i^{\cs}(p)$, so $\mse p\in Q$ and hence $\mse p\in W$.
\end{proof}
Thus it only remains to show $\mse n_1$ and $\mse m_n$ satisfy the
iterability conditions.  We will give the proof for $\mse n_1$; the
proof for $\mse m_n$ is similar.
\begin{claim}
\label{c.3}
The phalanx $(\tilde\gd,(W,\mse n_1))$ is iterable.
\end{claim}
\begin{proof}
We will suppose that $(\tilde\gd,(W,\mse n_1))$ is not iterable
and find an ill behaved tree on the phalanx $(\tilde\gd,(W,\mse n))$.
This contradicts the fact that  $(\tilde\gd,(W,\mse  n))$ is iterable,
since $\mse n$ is  a mouse of $W$, and thus proves claim~\ref{c.3}.  

The construction is similar to that in lemma~\ref{l.iterable2}.  
Let $\cs$ be an ill
behaved tree on  $(\tilde\gd,(W,\mse n_1))$, and let $Z$ be a
countable set containing everything relevant such that $Z\prec
H_{\tau}$ for some sufficiently large $\tau$. Let
$\eta=\sup(X\cap\tilde\gd)$ and let $k\colon \mse n'\to\mse n_1$ be
the inverse of the collapse map of the Skolem hull of
$\eta\cup\{\tilde\gd\}$ in  $\mse n_1$. 
As in lemma~\ref{l.iterable2} we can use the Martin-Steel iterability
theorem and a Levy collapse to show that the $\cs$ induces a ill
behaved tree $\cs'$ on $(\eta,(W,\mse n'))$.
Since $\card
Z=\gw<\gd=\cof(\tilde\gd)$ there is $\nu\in\pi\image I\setminus\eta$
such that $\nu\in_0 x$ for all $x\in Z\cap U$.  If $o(\nu)=0$ then we
can embed $\mse n'$ into $\mse n$ by mapping $\tilde\gd$ to $\nu$, and
thus we get an ill behaved tree on $(\eta,(W,\mse n'))$, a
contradiction. 

If $o(\nu)>0$ then we similarly find an ill behaved tree on
$(\eta,(W,\ult(\mse n,U_\nu)))$ where $U_\nu$ is the order~0
measure on $\nu$. This  finishes the proof of claim~\ref{c.3}.
\end{proof}
This completes the proof of claim~\ref{c.1} and hence of
lemma~\ref{UinW}. 
\end{proof}

We can now complete the proof of the main theorem.  It only remains to
use the ultrafilter $U$ on $m^*$ in $K$ to define a set of
indiscernibles for 
$\ca$ in $K$.
Let $h^*$ be the Skolem function  of
$\mse m^*$ and define $D$ to be  the set of ordinals $\ga\in C\cap\tilde\gd$
such that
\begin{gather*}
\forall x\in
h^*\image\ga\cap \ps(\tilde\gd)\;\left(\ga\in x\iff x\in U\right)\\
\intertext{and for all $n>0$ in $\gw$}
\forall x\in
h^*\image\ga\cap \ps([\tilde\gd]^{1+n})\;\left(
\set{\vec\nu\in[\tilde\gd]^n:(\ga,\vec\nu)\in x}\in U^{n}\iff
x\in U^{1+n}\right).
\end{gather*}
Then $D\in K$, and 
$\pi(\nu)\in_0 D$ for every sufficiently large ordinal
$\nu\in I$.
Thus every member $\nu$ of $\pi\image I$ either is in $D$ or has
$D\cap\nu\in U_{\nu}$, so $\card{D}\ge\card{\pi\image I}=\gd$.

This completes the proof of the main theorem.
\section{Some questions}

The most basic open question is the problem with which we opened the paper:
\begin{question}
Is is a theorem of ZFC that every Ramsey cardinal is \jonson?
\end{question}
It would certainly be surprising if this were true, but it is also
surprising that no counterexamples are yet known.

It is perhaps more plausible to hope that the restriction of the main
theorem to models with no class model of a Woodin cardinal can be
eliminated:
\begin{question}
Is it a theorem of $ZF$ that if the  core model $K$ exists then every
\jonson{} cardinal is Ramsey in $K$?
\end{question}
Of course the general notion of ``$K$ is the core model'' remains to
be defined.  For the present we could take the problem as refering to
core models in the sense of Steel.

Finally we mention one more question:
\begin{question}
Can theorem~\ref{t:j->erdos} be generalized to singular cardinals~$\gd$?
\end{question}
It is easy to see that some such generalization is possible, but it is
not clear how much can be said.


\end{document}